\documentclass{amsproc}
\usepackage{amssymb}
\usepackage{graphicx}
\usepackage{amscd}
\usepackage{amsmath}

\theoremstyle{plain}

\newtheorem{algorithm}{Algorithm}

\newtheorem{lemma}{Lemma}

\newtheorem{remark}{Remark}

\newtheorem{theorem}{Theorem}
\numberwithin{equation}{section}

\input{tcilatex}

\begin{document}
\title[]{The boundary control approach to the Titchmarsh-Weyl $m-$function.
I. The response operator and the $A-$amplitude }
\date{July, 2006}
\author{Sergei Avdonin, Victor Mikhaylov and Alexei Rybkin}
\address{Department of Mathematics and Statistics \\
University of Alaska Fairbanks\\
PO Box 756660\\
Fairbanks, AK 99775}
\email{ffsaa@uaf.edu, ftvsm@uaf,edu, ffavr@uaf.edu}
\maketitle

\begin{abstract}
We link the Boundary Control Theory and the Titchmarsh-Weyl Theory. This
provides a natural interpretation of the $A-$amplitude due to Simon and
yields a new efficient method to evaluate the Titchmarsh-Weyl $m-$function
associated with the Schr\"{o}dinger operator $H=-\partial _{x}^{2}+q\left(
x\right) $ on $L_{2}\left( 0,\infty \right) $ with Dirichlet boundary
condition at $x=0.$
\end{abstract}

\begin{center}
{\small \textit{Dedicated with great pleasure to B.S. Pavlov on the occasion
of his 70th birthday}}
\end{center}

\section{\protect\bigskip \textbf{Introduction }}

Consider the Schr\"odinger operator

\begin{equation}
H=-\partial _{x}^{2}+q\left( x\right)  \label{1.1}
\end{equation}
on $L_{2}\left( \mathbb{R}_{+}\right) \,,\mathbb{R}_{+}:=[0,\infty ),$ with
a real-valued locally integrable potential $q.$ We assume that $\left( \ref%
{1.1}\right) $ is limit point case at $\infty $, that is, for each $z\in
\mathbb{C}_{+}:=\{z\in \mathbb{C}:\func{Im}z>0\}$ the equation
\begin{equation}
-u^{\prime \prime }+q\left( x\right) u=zu  \label{1.2}
\end{equation}
has a unique, up to a multiplicative constant, solution $u_{+}$ which is in $%
L_{2}$ at $\infty ~$:
\begin{equation}
\int_{\mathbb{R}_{+}}\left\vert u_{+}\left( x,z\right) \right\vert
^{2}dx<\infty ,z\in \mathbb{C}_{+}.  \label{1.3}
\end{equation}
Such solution $u_{+}$ is called a \emph{Weyl} \emph{solution} and its
existence for a very broad class of real potentials $q$ is the central point
of the Titchmarsh-Weyl theory.

The (principal or Dirichlet) \emph{Titchmarsh-Weyl }$m$\emph{-function}, $%
m(z),$ is defined for $z\in \mathbb{C}_{+}$ as
\begin{equation}
m\left( z\right) :=\frac{u_{+}^{\prime }\left( 0,z\right) }{u_{+}\left(
0,z\right) }.  \label{1.4}
\end{equation}%
The function $m(z)$ is analytic in $\mathbb{C}_{+}$ and satisfies the
Herglotz property:
\begin{equation}
m:\mathbb{C}_{+}\rightarrow \mathbb{C}_{+},  \label{1.5}
\end{equation}%
so $m$ satisfies a Herglotz representation theorem,
\begin{equation}
m\left( z\right) =c+\int_{\mathbb{R}}\left( \frac{1}{t-z}-\frac{t}{1+t^{2}}%
\right) d\mu \left( t\right) ,  \label{1.6}
\end{equation}%
where $c=\func{Re}m\left( i\right) ~$and $\mu $ is a positive measure
subject to
\begin{equation}
\int_{\mathbb{R}}\frac{d\mu \left( t\right) }{1+t^{2}}<\infty ,  \label{1.7}
\end{equation}%
\begin{equation}
d\mu \left( t\right) =\text{w-}\lim_{\varepsilon \rightarrow +0}\frac{1}{\pi
}\func{Im}m\left( t+i\varepsilon \right) dt.  \label{1.8}
\end{equation}

It is a fundamental fact of the spectral theory of ordinary differential
operators that the measure $\mu $ is the spectral measure of the Schr\"{o}%
dinger operator $\left( \ref{1.1}\right) $ with a Dirichlet boundary
condition at $x=0.$ Another fundamental fact is the \emph{Borg-Marchenko
uniqueness theorem} \cite{Borg52}, \cite{Marchenko50} stating
\begin{equation}
m_{1}=m_{2}\Longrightarrow q_{1}=q_{2}.  \label{1.9}
\end{equation}%
There is no explicit formula realizing $\left( \ref{1.9}\right) $ but there
are some Gelfand-Levitan-Marchenko type procedures to recover the potential $%
q$ by given $m-$function.

The Titchmarsh-Weyl $m-$function is a central object of the spectral theory
of linear ordinary differential operators but its actual computation is
problematic. In fact, $\left( \ref{1.6}\right) $ is suitable if the spectral
measure $\mu $ of $\left( \ref{1.1}\right) $ with Dirichlet boundary
condition at $0$ is available, which is not usually the case. Instead, $%
\left( \ref{1.6}\right) $ is used to find $\mu $ by $\left( \ref{1.8}\right)
$ but not the other way around.

The definition $\left( \ref{1.4}\right) $ is not always practical either
since finding $m\left( z\right) $ by $\left( \ref{1.4}\right) $ is
essentially equivalent to solving
\begin{equation}
\left\{
\begin{array}{c}
-u^{\prime \prime }+q\left( x\right) u=zu \\
u\left( 0,z\right) =1,\int_{0}^{\infty }\left\vert u\left( x,z\right)
\right\vert ^{2}dx<\infty
\end{array}%
\right.   \label{1.9'}
\end{equation}%
for all $z\in \mathbb{C}_{+}.$ The analysis of the asymptotic behavior of $%
m\left( z\right) $ for large $\left\vert z\right\vert $ has received
enormous attention and the picture is now quite clear (see, e.g. \cite{CG'02}%
, \cite{Ry'02} and the literature cited therein). If $x=0$ is a (right)
Lebesgue point of $q\left( x\right) $ then%
\begin{equation}
m\left( z\right) =i\sqrt{z}+\frac{q\left( 0\right) }{2i\sqrt{z}}+o\left(
\frac{1}{\sqrt{z}}\right) ,~z\rightarrow \infty ,~\varepsilon \leq \arg
z\leq \pi -\varepsilon ,~\varepsilon >0,  \label{1.10}
\end{equation}%
which means that the $m-$functions for all $q$ coinciding on $\left[ 0,a%
\right] $ with arbitrarily small $a>0$ have the same asymptotic behavior.
Due to $\left( \ref{1.9}\right) $, it is therefore $m\left( z\right) $ for
finite $z$ that is of particular interest, which requires a very accurate
control of the solution to $\left( \ref{1.9'}\right) $ at $x\rightarrow
\infty $. In other words, the main issue here is the asymptotic behavior of $%
u\left( x,z\right) $ as $x\rightarrow \infty $ for finite $z.$ Typically,
such asymptotics are derived by transforming $\left( \ref{1.9'}\right) $ to
a suitable linear Volterra type integral equation. This can efficiently be
done when, e.g., $q$ decays at $\infty $ fast enough ($q\in L_{1}\left(
\mathbb{R}_{+}\text{ }\right) $ is sufficient). Equation $\left( \ref{1.9'}%
\right) $ can then be transformed to
\begin{equation*}
y\left( x,z\right) =1+\int_{x}^{\infty }K\left( x,s,z\right) y\left(
s,z\right) ds,~y\left( x,z\right) :=e^{-i\sqrt{z}x}u\left( x,z\right) ,
\end{equation*}%
where
\begin{equation*}
K\left( x,s,z\right) :=\frac{e^{-2i\sqrt{z}\left( s-x\right) }-1}{2i\sqrt{z}}%
q\left( x\right) ,
\end{equation*}%
which can be solved by iteration.

Another well-known transformation of $\left( \ref{1.2}\right) $ is the
Green-Liouville transformation $~$(see, e.g. \cite{T})
\begin{equation}
y^{\prime \prime }\left( x,z\right) +y\left( x,z\right) +\left[ \frac{1}{4}%
\frac{q^{\prime \prime }\left( x\right) }{\left\{ z-q\left( x\right)
\right\} ^{2}}+\frac{5}{16}\frac{q^{\prime 2}\left( x\right) }{\left\{
z-q\left( x\right) \right\} ^{3}}\right] y\left( x,z\right) =0,  \label{1.11}
\end{equation}%
where $y\left( x,z\right) =\left\{ z-q\left( x\right) \right\} ^{1/4}u\left(
x,z\right) $. Equation $\left( \ref{1.11}\right) $ is a crucial ingredient
in the WKB-analysis and can be reduced to a linear Volterra integral
equation for a wide range of potentials (even growing at infinity) but
requires that $q$ be twice differentiable. Even for smooth potentials like $%
q\left( x\right) =x^{-\alpha }\sin x^{\beta },0<\alpha \leq \beta \leq 1,$
the transformation $\left( \ref{1.11}\right) $ is not of much help since $%
q^{\prime }$ and $q^{\prime \prime }$ unboundedly oscillate at $\infty .$
Note, that, as it was shown by Buslaev-Matveev \cite{BM70}, the
Green-Liouville transformation $\left( \ref{1.11}\right) $ works well for
slowly decaying potentials subject to
\begin{equation}
\left\vert q^{\left( l\right) }\left( x\right) \right\vert \leq Cx^{-\alpha
-l},~\alpha >0,~l=0,1,2.  \label{1.12}
\end{equation}

One of the authors \cite{Ry'04} has recently put forward yet another
transformation that allows one to obtain and analyze the asymptotics for the
solution to $\left( \ref{1.9'}\right) $ for general non-smooth potentials $q$
with a very mild decay at $\infty $. Namely, if a potential $q$ is such that
the sequence $\left\{ \int_{j}^{j+1}\left\vert q\left( x\right) \right\vert
dx\right\} _{j=1}^{\infty }\,\ $is from $l_{p},p=2^{n}$ with some $n\in
\mathbb{N}$ then $\left( \ref{1.9'}\right) $ can be transformed to
\begin{equation}
y\left( x,z\right) =1+\int_{x}^{\infty }K_{n}\left( x,s,z\right) y\left(
s,z\right) ds,  \label{1.13}
\end{equation}%
where $y\left( x,\lambda \right) :=\Theta _{n}^{-1}\left( x,\lambda \right)
u\left( x,\lambda \right) $ and
$$
\Theta _{n}\left( s,z\right):=\Theta _{n}\left( 0,s,z\right)
,\Theta
_{n}\left( x,s,z\right):=\exp \left\{ i\sqrt{z}s+\int_{x}^{s+x}%
\sum_{m=1}^{n}q_{m}\left( t,z\right) dt\right\} ,$$
$$
K_{n}\left( x,s,\lambda \right):=\left( q_{n}\Theta _{n}\right)
^{2}\left( s,\lambda \right) \int_{x}^{s}\Theta _{n}^{-2}\left(
t,\lambda \right) dt.$$ The functions $q_{m}$ are, in turn,
defined by the following recursion formulas:
\begin{eqnarray}
q_{1}\left( x,z\right)  &:&=-\int_{0}^{\infty }e^{2i\sqrt{z}s}q\left(
s+x\right) ds,  \label{1.14} \\
~q_{m+1}\left( x,z\right)  &:&=\int_{0}^{\infty }\Theta _{m}^{2}\left(
x,s,z\right) q_{m}^{2}\left( s+x,z\right) ds,~m\in \mathbb{N}.  \notag
\end{eqnarray}%
Formulas $\left( \ref{1.14}\right) $ can be viewed as \textquotedblleft
energy dependent\textquotedblright\ transformations of the original
potential $q$ improving its rate of decay at infinity. For $n\geq 2$ these
transformations are highly nonlinear and were previously considered by many
authors (see, e.g. \cite{H2}, \cite{H}, \cite{K1}) in connection with a
variety of improvements of asymptotics $\left( \ref{1.10}\right) $. The main
feature of the transformations $\left( \ref{1.13}\right) -\left( \ref{1.14}%
\right) $ is that they yield higher order WKB type asymptotics of
the Weyl solution $u_{+}\left( x,z\right) $ as $x\rightarrow
\infty $ at fixed finite $z$ \ (see \cite{Ry'04}).

Our list of transformations of the original equation $\left( \ref{1.2}%
\right) $ is of course incomplete and given here to demonstrate how
drastically computational complexity\ of solving $\left( \ref{1.9'}\right) $
(and hence the $m-$function) tends to increase when one relaxes decay
conditions. It should be particularly emphasized that in order to get $%
m\left( z\right) $ one has to solve the integral equations for
each $z$ separately. In addition, $q$'s with no decay at $\infty $
should be considered on {\it ad hoc} basis.

In the present note we put forward a different approach to evaluate the
Titchmarsh-Weyl $m-$function which is based on the Boundary Control Theory.
The main idea is that to study the (dynamic) Dirichlet-to-Neumann map $%
u\left( 0,t\right) \mapsto u_{x}\left( 0,t\right) $ for the wave equation
associated with $\left( \ref{1.2}\right) $:
\begin{equation*}
u_{tt}-u_{xx}+q\left( x\right) u=0,\ x>0,~t>0,
\end{equation*}%
with zero initial conditions. The Dirichlet-to-Neumann map defined this way
turns out to be the so-called response operator, an important object of the
Boundary Control method in inverse problems \cite{ABI91}, \cite{AB96}. In
the frequency domain the latter becomes the operator of multiplication by
the Titchmarsh-Weyl $m$-function\footnote{%
Usually referred to as the Dirichlet (or principal) Titchmarsh-Weyl $m-$%
function} associated with the operator $-\partial _{x}^{2}+q\left( x\right) $
with Dirichlet boundary condition at $x=0$ (Pavlov \cite{Pavlov} has noticed
that $m$-function can be interpreted as a one-dimensional (spectral)
Dirichlet-to-Neumann map). This approach allows one to employ powerful
techniques developed for the wave equation to the study of the
Titchmarsh-Weyl $m-$function. In this paper we concentrate on the direct
problem only. That is, given potential $q$, we evaluate the $m-$function in
terms of the response operator (response function, to be precise) which is
exactly Simon's representation of the $m-$function via his $A-$amplitude
(see \cite{S'99} and \cite{GS'00}). Our approach however provides a clear
physical interpretation of the $A-$amplitude and gives a new procedure to
compute it. The latter can potentially be used for numerical analysis of the
$m-$function.

We emphasize that all the ingredients we use in the present paper are
already known in different inverse problems communities (\cite{GS'00}, \cite%
{LeSa'75}, \cite{KKLM}, \cite{SaCh'89}, \cite{FY'01} to name just five) but
it is the new way to combine them that makes our main contribution to this
well developed area. However, we do not utilize here the full power of the
Boundary Control approach, which is in inverse methods. We plan to address
this important issue in our sequel on this topic.

The paper is organized as follows. In Section 2 we introduce the main
ingredient of our approach, the response operator $R$, and give its
connection with the Titchmarsh-Weyl $m-$function. We also show that its
kernel is closely related to the $A-$amplitude.

In Section 3 we derive a linear Volterra type integral equation for a
function $A\left( x,y\right) $ which diagonal value is the $A-$amplitude
(Theorem 1).

Section 4 is devoted to the analysis of the integral equation for the kernel
$A\left( x,y\right) $ producing an important bound for the $A-$amplitude
(Theorem 2) which answers an open question by Gesztesy-Simon \cite{GS'00}.

In short Section 5 we present our algorithm of practical evaluation of the $%
m-$function and make some concluding remarks.

\section{The response operator and the $A$-amplitude}

Let us associate with the Schr\"{o}dinger equation $(\ref{1.1})$ the
axillary wave equation
\begin{equation}
\left\{
\begin{array}{c}
u_{tt}(x,t)-u_{xx}(x,t)+q(x)u(x,t)=0,\quad x>0,\ t>0, \\
u(x,0)=u_{t}(x,0)=0,\ u(0,t)=f(t),%
\end{array}%
\right.   \label{2.1'}
\end{equation}%
where $f$ is an arbitrary $L_{2}\left( \mathbb{R}_{+}\right) $ function
referred to as a \emph{boundary control}. It can be verified by a direct
computation that the weak solution $u^{f}(x,t)$ to the initial-boundary
value problem $(\ref{2.1'})$ admits the representation
\begin{equation}
u^{f}(x,t)=\left\{
\begin{array}{c}
f(t-x)+\int_{x}^{t}w(x,s)f(t-s)\,ds,~x\leq t, \\
0,\quad x>t,%
\end{array}%
\right.   \label{2.2}
\end{equation}%
in terms of the solution $w(x,s)$ to the Goursat problem:
\begin{equation}
\left\{
\begin{array}{c}
w_{ss}(x,s)-w_{xx}(x,s)+q(x)w(x,s)=0,\quad 0<x<s, \\
w(x,0)=0,w(x,x)=-\frac{1}{2}\int_{0}^{x}q\left( t\right) dt.%
\end{array}%
\right.   \label{2.3}
\end{equation}%
We introduce now the \emph{response operator} $R$:
\begin{equation}
(Rf)(t)=u_{x}(0,t),  \label{2.3'}
\end{equation}%
so it transforms $u\left( 0,t\right) \mapsto u_{x}\left(
0,t\right) $. By this reason it can also be called the (dynamic)
Dirichlet-to-Neumann map. From $(\ref{2.2})$ we easily get the
representation
\begin{equation}
(Rf)(t)=-\frac{d}{dt}f(t)+\int_{0}^{t}r(t-s)f(s)\,ds\,,
\label{2.4}
\end{equation}%
\begin{equation}
r(\cdot ):=w_{x}(0,\cdot ).  \label{2.4'}
\end{equation}%
In other words, the response operator is the operator of differentiation
plus the convolution. The kernel $r$ of the convolution part of $\left( \ref%
{2.4}\right) $ is called the \emph{response function} which plays
an important role in the Boundary Control method (see, e.g.
\cite{ABI91}, \cite{AB96}).

In fact, besides the Boundary Control method some close analogs of
the response function have independently been discovered in the
half-line short-range scattering \cite{Ramm'87} and more recently
in the connection
with inverse spectral problem for the half-line Schrodinger operator \cite%
{S'99}, \cite{GS'00}, \cite{RaSi'00}.

We now demonstrate the connection between the response function $r\left(
s\right) $ and the (Dirichlet) Titchmarsh-Weyl $m$-function. An interplay
between spectral and time-domain data  is widely used in inverse problems,
see, e.g., \cite{KKLM} where  the equivalence of several types of boundary
inverse problems is  discussed for smooth coefficients; notice, however,
that we consider the case of not smooth but  just $L^1_{loc}$ potentials.

Let $f \in C_0^{\infty}(0,\infty)$ and
\begin{equation*}
\widehat{f}(k):=\int_0^{\infty} f(t)\,e^{-kt}\,dt\,
\end{equation*}
be its Laplace transform. Function $\widehat{f}(k)$ is well defined for $k
\in \mathbb{C}$ and, if $\func{Re} k >0\,,$
\begin{equation}  \label{est}
|\widehat{f}(k)| \leq C_\alpha (1+|k|)^{-\alpha}
\end{equation}
for any $\alpha>0.$ Going in $(\ref{2.1'})\,$ and $(\ref{2.3'})\,$ over to
the Laplace transforms, one has
\begin{eqnarray}
-\widehat{u}_{xx}(x,k)+q(x)\widehat{u}(x,k) &=&-k^{2}\widehat{u}(x,k),
\label{2.6} \\
\widehat{u}(0,k) &=&\widehat{f}(k),  \label{2.7}
\end{eqnarray}%
and
\begin{equation}
\widehat{(Rf)}(k)=\widehat{u}_{x}(0,k)\,,
\end{equation}%
respectively.

Estimate (\ref{est}) implies that $|\widehat{u}(x,k)|$ decreases rapidly
when $|k| \rightarrow \infty\,,\; \func{Re} k \geq \epsilon >0\,.$ The
values of the function $\widehat{u} (0,k) $ and its first derivative at the
origin, $\widehat{u}_x (0,k), $ are related through the Titchmarsh-Weyl
m-function
\begin{equation}
\widehat{u}_x (0,k)=m(-k^2)\widehat{f}(k)\,.  \label{psim}
\end{equation}
Therefore,
\begin{equation}
\widehat{(Rf)}(k)=m(-k^2)\widehat{f}(k)\,,  \label{2.7'}
\end{equation}
and thus the spectral and dynamic Dirichlet-to-Neumann maps are in
one-to-one correspondence.

Taking the Laplace transform of $(\ref{2.4})\,$\ we get
\begin{equation}
\widehat{(Rf)}(k) =-k\widehat{f}(k)+\widehat{r}(k)\widehat{f}(k).
\label{2.8}
\end{equation}%
In Section 4 we show that, under some mild conditions on the potential $q,$ $%
(\ref{2.7'})\,$\ and $(\ref{2.8})\,$\ imply
\begin{equation}
m(-k^{2})=-k+\int_{0}^{\infty }e^{-k\alpha }r(\alpha )\,d\alpha \, ,
\label{2.5}
\end{equation}%
where the integral is absolutely convergent in a proper domain of $k$.

Representation $(\ref{2.5})$ is not new. In the form
\begin{equation}
m(-k^{2})=-k-\int_{0}^{\infty }A(\alpha )e^{-2\alpha k}\,d\alpha  \label{2.9}
\end{equation}%
(with the absolutely convergent integral) it was proven for $q\in
L^{1}\left( \mathbb{R}_{+}\right) $ and $q\in L^{\infty }\left( \mathbb{R}%
_{+}\right) $ by Gesztesy-Simon \cite{GS'00} who call the function $A$ in\ $(%
\ref{2.9})$ the\emph{\ }$A-$\emph{amplitude}. Clearly, one has
\begin{equation}
A(\alpha )=-2r(2\alpha ).  \label{2.10}
\end{equation}

\begin{remark}
\bigskip The fundamental role of representation $(\ref{2.9})$ and the $A-$%
amplitude was emphasized in \cite{S'99} and \cite{GS'00}. However no direct
interpretation of $(\ref{2.9})$ and $A$ is given in \cite{S'99}, \cite{GS'00}%
. On the other hand, $(\ref{2.5})$ says that the Titchmarsh-Weyl $m-$%
function is the Laplace transform of the kernel of the response operator $%
R\,\ $(see $\left( \ref{2.4}\right) $). Or, equivalently, the
matrix (one by one in our case) of the response operator $R$ in
the spectral representation of $H_{0}=-\partial _{x}^{2},u\left(
0\right) =0,$ coincides with the Titchmarsh-Weyl $m-$function
associated with $H=H_{0}+q.$ The response operator $R,$ in turn,
describes the reaction of the system. In particular, for the
string the operator $R$ connects the displacement and tension at
the endpoint $x=0.$ For electric circuits it relates the current
and voltage (see, e.g. \cite{Pavlov}, \cite{KKLM} and references
therein for additional information about the physical meaning of
the Dirichlet-to-Neumann map). In the theory of linear dynamical
systems the response operator is the input-output map and the
Laplace transform of its kernel is the transfer function of a
system.
\end{remark}

\begin{remark}
In fact, $\left( \ref{2.7'}\right)\, $ can be regarded as a
definition of the Titchmarsh-Weyl $m-$function which could be
effortlessly extended to matrix valued and complex potentials
since the Boundary Control method is readily available in these
situations {\cite{ABI91, AB96}. We hope to return to this
important point elsewhere. Also, since the Dirichlet-to-Neumann
map can be viewed as a $3D$ analog of the $m-$function, $\left(
\ref{2.7'}\right) $ could hopefully yield a canonical way to
define (operator valued) $m-$functions for certain partial
differential operators. It is worth mentioning that Amrein-Pearson
\cite{AmPe'04} have recently generalized (using quite different
methods) the theory of the Weyl-Titchmarsh $m-$function for
second-order ordinary differential operators to partial
differential operators of the form $-\Delta +q\left( x\right) $
acting in three space dimensions.$\ $ }
\end{remark}

Despite the clear physical interpretation of the response function $r$ some
formulas in Section 3 look slightly prettier in terms of the $A-$amplitude.
Since our interest to the topic was originally influenced by \cite{S'99},
\cite{GS'00} we therefore are going to deal with $A$ related to $r$ by $%
\left( \ref{2.10}\right) $. In Section 4 we prove the absolute convergence
of the integral in $(\ref{2.9})$ (and, therefore, of the integral in $(\ref%
{2.5})$) for
\begin{equation*}
q\in l^{\infty }\left( L^{1}\left( \mathbb{R}_{+}\right) \right) :=\left\{
q:\int_{n}^{n+1}\left\vert q\left( x\right) \right\vert \,dx\in l^{\infty
}\right\} .
\end{equation*}

\section{An integral equation for the $A$-amplitude}

In this section we derive a linear Volterra type integral equation closely
related to the $A-$amplitude.

\begin{theorem}
Let $q\in L_{\limfunc{loc}}^{1}\left( \mathbb{R}_{+}\right) .$ Then for a.
e. $\alpha >0$
\begin{equation}
A\left( \alpha \right) =A\left( \alpha ,\alpha \right) ,  \label{3.1}
\end{equation}%
where $A\left( x,y\right) $ is the solution to the integral equation
\begin{equation}
A\left( x,y\right) =q\left( x\right) -\int_{0}^{y}\left(
\int_{v}^{x}\,A\left( u,v\right) du\right) q\left( x-v\right) dv \,; \
x,y>0\,.  \label{3.2}
\end{equation}
\end{theorem}

\begin{proof}
We go through a chain of standard transformations of the Goursat problem $%
\left( \ref{2.3}\right) $. By setting $u=s+x,v=s-x$ and
\begin{equation}
V\left( u,v\right) =w\left(\frac{u-v}{2},\frac{u+v}{2}\right)\,,
\label{3.2'}
\end{equation}
equation $\left( \ref{2.3}\right) $ reduces to
\begin{equation*}
\left\{
\begin{array}{c}
V_{uv}+4q\left( \dfrac{u-v}{2}\right) V=0 \\
V\left( u,u\right) =0 \\
V\left( u,0\right) =-\frac{1}{2}\int_{0}^{u/2}q%
\end{array}
\right.
\end{equation*}
which can be easily transformed into
\begin{equation}
V\left( u,v\right) =-\frac{1}{2}\int_{v/2}^{u/2}q-\frac{1}{4}%
\int_{0}^{v}dv_{1}\int_{v}^{u}du_{1}q\left( \dfrac{u_{1}-v_{1}}{2}\right)
V\left( u_{1},v_{1}\right) .  \label{3.3}
\end{equation}
Doubling the variables in $\left( \ref{3.3}\right) $ yields
\begin{equation}
V\left( 2u,2v\right) =-\frac{1}{2}\int_{v}^{u}q-\int_{0}^{v}dv_{1}%
\int_{v}^{u}du_{1}q\left( u_{1}-v_{1}\right) V\left( 2u_{1},2v_{1}\right) .
\label{3.4}
\end{equation}
Introduce a new function
\begin{equation}
U\left( x,y\right) :=\int_{0}^{y}dvq\left( x-v\right) V\left( 2x,2v\right) .
\label{3.4'}
\end{equation}
It follows from $\left( \ref{3.4}\right) $ that $U\left( x,y\right) $
satisfies the integral equation
\begin{equation*}
U\left( x,y\right) =-\frac{1}{2}\int_{0}^{y}dv\,q\left( x-v\right)
\int_{v}^{x}du\,q\left( u\right) -\int_{0}^{y}dv\,q\left( x-v\right)
\int_{v}^{x}du\,U\left( u,v\right) .
\end{equation*}
The function
\begin{equation}
A\left( x,y\right) =q\left( x\right) +2U\left( x,y\right)  \label{3.4''}
\end{equation}
then obeys equation $\left( \ref{3.2}\right) .$

It is left to show $\left( \ref{3.1}\right) .$ By $\left( \ref{2.10}\right) $
and $\left( \ref{2.4'}\right) $
\begin{equation}
A\left( \alpha \right) =-2r\left( 2\alpha \right) =-2w_{x}\left( 0,2\alpha
\right) .  \label{3.5}
\end{equation}
But it follows from $\left( \ref{3.2'}\right) $ that
\begin{equation*}
w_{x}\left( x,s\right) =\left( V_{u}-V_{v}\right) \left( s+x,s-x\right)
\end{equation*}
and hence
\begin{equation}
w_{x}\left( 0,2\alpha \right) =\left( V_{u}-V_{v}\right) \left( 2\alpha
,2\alpha \right) .  \label{3.6}
\end{equation}
Differentiating $\left( \ref{3.4}\right) $ with respect to $u$ and $v$ and
setting $u=v=2\alpha $, we have
\begin{eqnarray*}
V_{u}\left( 2\alpha ,2\alpha \right) &=&-\frac{1}{4}q\left( \alpha \right) -%
\frac{1}{4}\int_{0}^{2\alpha }dv_{1}q\left( \alpha -\frac{v_{1}}{2}\right)
V\left( 2\alpha ,v_{1}\right) , \\
V_{v}\left( 2\alpha ,2\alpha \right) &=&\frac{1}{4}q\left( \alpha \right) +%
\frac{1}{4}\int_{0}^{2\alpha }dv_{1}q\left( \alpha -\frac{v_{1}}{2}\right)
V\left( 2\alpha ,v_{1}\right) .
\end{eqnarray*}
Inserting these formulas into $\left( \ref{3.6}\right) $ we get
\begin{equation}
w_{x}\left( 0,2\alpha \right) =-\frac{1}{2}q\left( \alpha \right) -\frac{1}{2%
}\int_{0}^{2\alpha }dv_{1}q\left( \alpha -\frac{v_{1}}{2}\right) V\left(
2\alpha ,v_{1}\right) .  \label{3.7}
\end{equation}
Setting in $\left( \ref{3.7}\right) $ $v_{1}=2v$ and plugging it then in $%
\left( \ref{3.5}\right) $, yields
\begin{equation}
A\left( \alpha \right) =q\left( \alpha \right) +2\int_{0}^{\alpha
}dv\,q\left( \alpha -v\right) V\left( 2\alpha ,2v\right) .  \label{3.8}
\end{equation}
It is left to notice that by $\left( \ref{3.4'}\right) $ the right hand side
of $\left( \ref{3.8}\right) $ is $q\left( \alpha \right) +2U\left( 2\alpha
,2\alpha \right) $ which by $\left( \ref{3.4''}\right) $ is equal to $%
A\left( \alpha ,\alpha \right) $ and $\left( \ref{3.1}\right) $ \ is proven.
\end{proof}

The kernel $A\left( x,y\right) $ in Theorem 1 is not related to $A\left(
\alpha ,x\right) $ appearing in \cite{S'99}, \cite{GS'00} where $A\left(
\alpha ,x\right) $ is the $A-$amplitude corresponding to the $m-$function
associated with the interval $\left( x,\infty \right) $.

\section{\protect\bigskip Analysis of iterations}

In this section we demonstrate that integral equation $\left( \ref{3.2}%
\right) $ is quite easy to analyze. We need the following technical

\begin{lemma}
Let $f\left( x\right) $ be a non-negative function and
\begin{equation}
||f||:=\sup_{x\geq 0}\int_{x}^{x+1}f\left( s\right) ds<\infty .  \label{4.1}
\end{equation}%
Then for any $a,b\geq 0$ and natural $n$%
\begin{equation}
\int_{0}^{a}\left( x+b\right) ^{n}f\left( x\right) dx\leq \frac{\left(
a+b+1\right) ^{n+1}}{n+1}||f||.  \label{4.2}
\end{equation}
\end{lemma}

\begin{proof}
We may assume $||f||=1$. Integrating the left hand side of $\left( \ref{4.2}%
\right) $ by parts yields
\begin{eqnarray}
\int_{0}^{a}\left( x+b\right) ^{n}f\left( x\right) dx &=&-\int_{0}^{a}\left(
x+b\right) ^{n}d\left( \int_{x}^{a}f\left( s\right) ds\right)   \notag \\
&=&b^{n}\int_{0}^{a}f+\int_{0}^{a}\left( \int_{x}^{a}f\left( s\right)
ds\right) d\left( x+b\right) ^{n}.  \label{4.3}
\end{eqnarray}%
Due to the trivial inequality
\begin{equation*}
\int_{\alpha }^{\beta }f<\beta -\alpha +1,~
\end{equation*}%
$\left( \ref{4.3}\right) $\ can be estimated above as follows
\begin{eqnarray*}
\int_{0}^{a}\left( x+b\right) ^{n}f\left( x\right) dx &<&b^{n}\left(
a+1\right) +\int_{0}^{a}\left( a-x+1\right) d\left( x+b\right) ^{n} \\
&=&b^{n}\left( a+1\right) +\left( a+b\right) ^{n}+\frac{(a+b)^{n+1}}{n+1}%
-b^{n}\left( a+1\right) -\frac{b^{n+1}}{n+1} \\
&\leq &\left( a+b\right) ^{n}+\frac{(a+b)^{n+1}}{n+1} \\
&=&\frac{1}{n+1}\left\{ (a+b)^{n+1}+\left( n+1\right) \left( a+b\right)
^{n}\right\}  \\
&\leq &\frac{(a+b+1)^{n+1}}{n+1}.
\end{eqnarray*}%
At the last step we used the obvious inequality $\left( x\geq 0\right) $%
\begin{equation*}
\left( x+1\right) ^{n+1}\geq x^{n+1}+\left( n+1\right) x^{n}.
\end{equation*}
\end{proof}

The following theorem is the main result of this section.

\begin{theorem}
Let $q$ be subject to
\begin{equation}
||q||:=\sup_{x\geq 0}\int_{x}^{x+1}\left\vert q\left( s\right) \right\vert
ds<\infty \,.  \label{4.5}
\end{equation}%
Then for $\alpha \geq 0$
\begin{equation}
\left\vert A\left( \alpha \right) -q\left( \alpha \right) \right\vert \leq
\frac{1}{2}\left( \int_{0}^{\alpha }\left\vert q\left( s\right) \right\vert
ds\right) ^{2}\left\{ \exp \left( 2\sqrt{2}\sqrt{||q||}\alpha \right) +\frac{%
1}{\sqrt{2\pi }}\exp \left( 2e||q||\alpha \right) \right\} \,.  \label{4.6}
\end{equation}
\end{theorem}

\begin{proof}
Rewriting $\left( \ref{3.2}\right) $ as
\begin{equation*}
A=q-KA,
\end{equation*}%
where
\begin{equation*}
\left( Kf\right) \left( x,y\right) :=\int_{0}^{y}dv\,q\left( x-v\right)
\int_{v}^{x}du\,f\left( u,v\right) ,
\end{equation*}%
and formally solving it by iteration, we get
\begin{equation*}
A\left( \alpha \right) =q\left( \alpha \right) +\sum_{n\geq 1}\left(
-1\right) ^{n}A_{n}\left( \alpha \right) ,\ A_{n}\left( \alpha \right)
:=\left( K^{n}q\right) \left( \alpha ,\alpha \right) ,
\end{equation*}%
and hence
\begin{equation}
\left\vert A\left( \alpha \right) -q\left( \alpha \right)
\right\vert \leq \sum_{n\geq 1}\left\vert A_{n}\left( \alpha
\right) \right\vert \leq \sum_{n\geq 1}I_{n}\left( \alpha \right)
\,,~ \label{4.7}
\end{equation}%
where
\begin{equation}
I_{n}\left( \alpha \right) :=\left( \left\vert K\right\vert ^{n}\left\vert
q\right\vert \right) \left( \alpha ,\alpha \right)   \label{4.7'}
\end{equation}%
with the agreement that $\left\vert K\right\vert $ is the integral operator $%
K$ with $\left\vert q\right\vert $ in place of $q.$ We now need a suitable
estimate for $\left( \left\vert K\right\vert ^{n}\left\vert q\right\vert
\right) \left( x,y\right) $. For $n=1,$%
\begin{eqnarray*}
\left( \left\vert K\right\vert \left\vert q\right\vert \right) \left(
x,y\right)  &=&\int_{0}^{y}\left( \int_{v}^{x}\left\vert q\left( u\right)
\right\vert du\right) \left\vert q\left( x-v\right) \right\vert dv
 \\
&\leq &\int_{x-y}^{x}ds\left\vert q\left( s\right) \right\vert
\int_{0}^{x}dt\left\vert q\left( t\right) \right\vert \leq \left(
\int_{0}^{x}\left\vert q\left( s\right) \right\vert ds\right)
^{2}=:Q^{2}\left( x\right) .  \notag
\end{eqnarray*}%
Similarly,
\begin{eqnarray}
\left( \left\vert K\right\vert ^{2}\left\vert q\right\vert \right) \left(
x,y\right)  &=&\int_{0}^{y}\left( \int_{v}^{x}\left( \left\vert K\right\vert
\left\vert q\right\vert \right) \left( u,v\right) du\right) \left\vert
q\left( x-v\right) \right\vert dv  \notag \\
&\leq &\sup_{0\leq v\leq u\leq x}\left( \left\vert K\right\vert \left\vert
q\right\vert \right) \left( u,v\right) \,\int_{0}^{y}\left(
\int_{v}^{x}du\right) \left\vert q\left( x-v\right) \right\vert dv
\label{4.9} \\
&\leq &Q^{2}\left( x\right) \,\int_{0}^{y}\left( x-v\right) \left\vert
q\left( x-v\right) \right\vert dv  \notag \\
&\leq &Q^{2}\left( x\right) \,x\int_{0}^{y}\left\vert q\left( x-v\right)
\right\vert dv=Q^{2}\left( x\right) \,x\int_{0}^{y}\left\vert q\left(
v+(x-y)\right) \right\vert dv  \notag \\
&\leq &Q^{2}\left( x\right) \,x\left( y+1\right) \,||q||.  \notag
\end{eqnarray}%
Here the supremum in $\left( \ref{4.9}\right) $ was estimated by $\left( \ref%
{4.2}\right) .$ We are now able to make the induction assumption
\begin{equation}
\left( \left\vert K\right\vert ^{n}\left\vert q\right\vert \right) \left(
x,y\right) \leq Q^{2}\left( x\right) \frac{x^{n-1}}{\left( n-1\right) !}%
\frac{\left( y+n-1\right) ^{n-1}}{\left( n-1\right) !}||q||^{n-1}.
\label{4.9'}
\end{equation}%
Then
\begin{eqnarray*}
\left( \left\vert K\right\vert ^{n+1}\left\vert q\right\vert \right) \left(
x,y\right)  &\leq &\int_{0}^{y}\left( \int_{v}^{x}\left( \left\vert
K\right\vert ^{n}\left\vert q\right\vert \right) \left( u,v\right) du\right)
\left\vert q\left( x-v\right) \right\vert dv \\
\,\, &\leq &\int_{0}^{y}\left( \int_{v}^{x}Q^{2}\left( u\right) \frac{u^{n-1}%
}{\left( n-1\right) !}\frac{\left( v+n-1\right) ^{n-1}}{\left( n-1\right) !}%
||q||^{n-1}du\right) \left\vert q\left( x-v\right) \right\vert dv \\
&\leq &\sup_{0\leq u\leq x}Q^{2}\left( u\right) \left( \int_{0}^{x}\frac{%
u^{n-1}}{\left( n-1\right) !}du\right) \left( \int_{0}^{y}\frac{\left(
v+n-1\right) ^{n-1}}{\left( n-1\right) !}\left\vert q\left( x-v\right)
\right\vert dv\right) ||q||^{n-1} \\
&\leq &Q^{2}\left( x\right) \frac{x^{n}}{n!}\frac{\left( y+n\right) ^{n}}{n!}%
||q||^{n}.
\end{eqnarray*}%
At the last step we estimated the second integral by Lemma 1. One now
concludes that $\left( \ref{4.9'}\right) $ holds for all natural $n.$

Combining $\left( \ref{4.7'}\right) $ and $\left( \ref{4.9'}\right) $, we
get
\begin{equation*}
I_{n}\left( \alpha \right) =\left( \left| K\right| ^{n}\left| q\right|
\right) \left( \alpha ,\alpha \right) \leq Q^{2}\left( \alpha \right) \frac{%
\alpha ^{n-1}}{\left( n-1\right) !}\frac{\left( \alpha +n-1\right) ^{n-1}}{%
\left( n-1\right) !}||q||^{n-1}
\end{equation*}
and hence for $\left(\ref{4.7}\right) $ we have
\begin{equation}
\left| A\left( \alpha \right) -q\left( \alpha \right) \right| \leq
Q^{2}\left( \alpha \right) \sum_{n\geq 0}\frac{\alpha ^{n}}{n!}\frac{\left(
\alpha +n\right) ^{n}}{n!}||q||^{n}.  \label{4.10}
\end{equation}
By the inequality $\left( a+b\right) ^{n}\leq 2^{n-1}\left(
a^{n}+b^{n}\right) $, estimate $\left( \ref{4.10}\right) $ continues
\begin{eqnarray}
\left| A\left( \alpha \right) -q\left( \alpha \right) \right| &\leq
&Q^{2}\left( \alpha \right) \left\{ \sum_{n\geq 0}2^{n-1}\left( \frac{\alpha
^{n}}{n!}\right) ^{2}||q||^{n}+\sum_{n\geq 1}2^{n-1}\frac{\alpha ^{n}}{n!}%
\frac{n^{n}}{n!}||q||^{n}\right\}  \notag \\
&=&\frac{1}{2}Q^{2}\left( \alpha \right) \left\{ \sum_{n\geq 0}\left( \frac{%
\sqrt{2||q||}\alpha }{n!}\right) ^{2n}+\sum_{n\geq 1}\frac{\left( 2\alpha
\right) ^{n}}{n!}\frac{n^{n}}{n!}||q||^{n}.\right\}  \label{4.11}
\end{eqnarray}
For the first series on the right hand side of $\left( \ref{4.11}\right) $
one has
\begin{equation}
\sum_{n\geq 0}\left( \frac{\sqrt{2||q||}\alpha }{n!}\right) ^{2n}\leq \exp
^{2}\left( \sqrt{2||q||}\alpha \right) =\exp \left( 2\sqrt{2||q||}\alpha
\right) .  \label{4.12}
\end{equation}
Evaluate the other one. It follows from the Stirling formula that
\begin{equation*}
n!\geq \sqrt{2\pi }\sqrt{n}\left( \frac{n}{e}\right) ^{n}
\end{equation*}
and hence
\begin{eqnarray}
\sum_{n\geq 1}\frac{\left( 2\alpha \right) ^{n}}{n!}\frac{n^{n}}{n!}%
||q||^{n} &\leq &\frac{1}{\sqrt{2\pi }}\sum_{n\geq 1}\frac{\left( 2\alpha
\right) ^{n}}{n!}\frac{e^{n}}{\sqrt{n}}||q||^{n}  \notag \\
&=&\frac{1}{\sqrt{2\pi }}\sum_{n\geq 1}\frac{\left( 2e||q||\alpha \right)
^{n}}{n!}<\frac{1}{\sqrt{2\pi }}\exp \left( 2e||q||\alpha \right) .
\label{4.13}
\end{eqnarray}
It follows now from $\left( \ref{4.11}\right) -\left( \ref{4.13}\right) $
that
\begin{equation*}
\left| A\left( \alpha \right) -q\left( \alpha \right) \right| <\frac{1}{2}%
Q^{2}\left( \alpha \right) \left\{ \exp \left( 2\sqrt{2||q||}\alpha \right) +%
\frac{1}{\sqrt{2\pi }}\exp \left( 2e||q||\alpha \right) \right\}
\end{equation*}
and $\left( \ref{4.6}\right) $ is proven.
\end{proof}

\begin{remark}
The exponential bounds $\left( \ref{4.6}\right) $ on $A\left( \alpha \right)
$ in Theorem 2 can be easily improved for the cases $q\in L^{1}\left(
0,\infty \right) $ and $q\in L^{\infty }\left( 0,\infty \right) $:
\begin{equation}
|A(\alpha )-q(\alpha )|\leqslant Q^{2}(\alpha )e^{\alpha Q(\alpha )},Q\left(
\alpha \right) =\int_{0}^{\alpha }\left\vert q\left( x\right) \right\vert dx,
\label{4.14}
\end{equation}%
\begin{equation}
|A(\alpha )-q(\alpha )|\leqslant ||q||_{\infty }\sum_{n\geq 1}\frac{%
||q||_{\infty }^{n}\alpha ^{2n}}{n!\left( n+1\right) !},||q||_{\infty
}:=\sup_{0\leq x\leq \infty }\left\vert q\left( x\right) \right\vert ,
\label{4.15}
\end{equation}%
respectively. Bounds $\left( \ref{4.14}\right) $ and $\left( \ref{4.15}%
\right) $ were found in \cite{GS'00}. In \cite{GS'00} (Section 10)
Gesztesy-Simon also conjectured that $A$ has exponential bound for
potentials obeying $\left( \ref{4.5}\right) $. Theorem 2 gives an
affirmative answer to their conjecture.
\end{remark}

\section{A new procedure for evaluating the Titchmarsh-Weyl $m-$function}

In this short section we present an algorithm to evaluate the $m-$function.

\begin{algorithm}
Given real valued potential $q$ subject to $\left\vert \left\vert
q\right\vert \right\vert =\sup \int_{x}^{x+1}\left\vert q\left( s\right)
\right\vert ds<\infty \,,$ the $m-$function can be computed as follows:

1. Solve integral equation $\left( \ref{3.2}\right) $ for $A\left(
x,y\right) $ and evaluate $A\left( \alpha \right) $ by $\left( \ref{3.1}%
\right) .$

2. Evaluate $m\left( z\right) $ by $\left( \ref{2.9}\right) $. The integral
in $\left( \ref{2.9}\right) $ is absolute convergent for $z=-k^{2}$ where $%
\func{Re}k>2\max \{\sqrt{2\left| \left| q\right| \right| },e\left| \left|
q\right| \right| \}.$
\end{algorithm}

\begin{remark}
Our algorithm yields an absolutely convergent series representation of the $%
m-$function,
\begin{equation}
m(-k^{2})=-k-\sum_{n\geq 0}\left( -1\right) ^{n}\int_{0}^{\infty
}A_{n}(\alpha )e^{-2\alpha k}\,d\alpha , \ A_{0}:=q.  \label{5.1}
\end{equation}%
Under weaker conditions on $q$ representation $\left( \ref{5.1}\right) $ was
obtained in \cite{GS'00} by completely different methods which do not imply
the linear Volterra integral equation $\left( \ref{3.2}\right) $. Some other
series representations can be found in \cite{H2}, \cite{H}, \cite{K1}, \cite%
{Ry'02}, \cite{Ry'04}. It should be pointed out though that those series are
quite unwieldy and, in addition, should be computed for each $z$ separately.
Our procedure has the advantage that once $A\left( \alpha \right) $ is found
one only needs to compute its Laplace transform for different $k$.
\end{remark}

\begin{remark}
It can be easily seen that if $q\left( x\right) \geq 0\,\ $ then $%
A_{n}(\alpha )\geq 0$ which improves the rate of convergence of $A\left(
\alpha \right) =\sum_{n\geq 0}\left( -1\right) ^{n}A_{n}(\alpha )$ making
our algorithm more efficient.
\end{remark}

\section{Acknowledgments}

\noindent Research of Sergei Avdonin was partially supported by the National
Science Foundation, grant \# OPP-0414128. Victor Mikhaylov was partially
supported by the University of Alaska Fairbanks Graduate Fellowship.

\end{document}